# PREHOMOGENEOUS VECTOR SPACES AND ERGODIC THEORY I

AKIHIKO YUKIE[1]

Oklahoma State University

**Contents**

**Introduction**

This is part one of a series of papers. In this series of papers, we consider problems analogous to the Oppenheim conjecture from the viewpoint of prehomogeneous vector spaces.

Throughout this paper, $k$ is a field of characteristic zero. The following theorem, known as the Oppenheim conjecture, was proved by Margulis [19].

**Theorem (0.1) (Margulis)** *Let $Q$ be a real non-degenerate indefinite quadratic form in $n \geq 3$ variables. Suppose that the corresponding point in $\mathbb{P}(\mathrm{Sym}^2(\mathbb{R}^n)^*)$ is irrational. Then the set of values of $Q$ at primitive integer points is dense in $\mathbb{R}$.*

The above theorem for $n \geq 5$ was conjectured by Oppenheim in [21]. Margulis originally proved that values of $Q$ at integer points can be arbitrarily small (non-trivially of course), which implies that the set of values of $Q$ at integer points is dense in $\mathbb{R}$ due to the result of Lewis [18]. The above improved version (the primitive part) is due to Dani–Margulis [8]. A further improvement of this result was obtained by Borel–Prasad [4]. Some partial results were known prior to the work of Margulis.

Let $f(x)$ be a degree $d$ form in real $n$ variables $x = (x_1, \cdots, x_n)$. In the following, we always assume that $f$ is not a multiple of an integral form. Consider the following questions.

(1) For any $\epsilon > 0$, does there exist $x \in \mathbb{Z}^n \setminus \{0\}$ such that $|f(x)| < \epsilon$?
(2) For any $\epsilon > 0$, does there exist $x \in \mathbb{Z}^n \setminus \{0\}$ such that $0 < |f(x)| < \epsilon$?

---

[1]Partially supported by NSF grant DMS-9401391

(3) Is the set $\{f(x) \mid x \in \mathbb{Z}^n\}$ dense in $\mathbb{R}$?

For non-degenerate quadratic forms, one typically gets stronger results such as (2) or (3). For forms of higher degree, there is no notion of "non-degenerate" forms. However, for (1), this is not necessary. For example, if the form does not depend on the variable $x_1$, we may choose $x_1 = 1$ and $x_2 = \cdots = x_n = 0$. For non-degenerate quadratic forms, the result of Lewis [18] implies that (2) implies (3).

Oppenheim himself had partial results for quadratic forms in [22], [23]. For example, he proved that if $f$ is a non-degenerate indefinite quadratic form in five or more variables and there exists $x \in \mathbb{Z}^n \setminus \{0\}$ such that $f(x) = 0$, then (2) is true. Chowla [5] proved (2) for indefinite diagonal quadratic forms in nine variables. Davenport and Heilbronn [12] proved (2) for indefinite diagonal quadratic forms in five variables. Davenport with the collaboration of Birch and Rigout proved (2) for indefinite quadratic forms in 21 or more variables (see [9], [10], [11], [2], [30]). For diagonal forms of degree $d$, Birch and Vaughan pointed out [3] that the method in [12] applies to forms in $2^d + 1$ variables. Davenport and Roth [13] proved a stronger result that if $d > 12$, $Cd \log d$ variables are enough ($C$ is a constant) and that if $d = 3$, eight variables are enough. Also for non-diagonal forms in higher degree, Schmidt [33] proved (1) for odd degree forms in enough variables. For example, for cubic forms, $(1314)^{256}$ variables are enough due to the result of Pitman [24]. Note that the result of Lewis is for quadratic forms, and for generic odd degree forms in sufficiently many variables, we still don't know whether or not (2) implies (3).

The method of Davenport and others are based on the circle method, and the method of Schmidt is somewhat more transcendental. However, the proof of Margulis is based on ergodic theory. What Margulis did was to prove the following theorem for the case $G = \mathrm{SL}(3)_\mathbb{R}$, $U = \mathrm{SO}(Q)_\mathbb{R}$ where $Q$ is a non-degenerate indefinite quadratic form in three variables.

**Theorem (0.2) (Ratner)** *Let $G$ be a connected Lie group and $U$ a connected subgroup of $G$ generated by unipotent elements of $G$. Then given any lattice $\Gamma \subset G$ and $x \in G/\Gamma$, there exists a connected closed subgroup $U \subset F \subset G$ such that $\overline{Ux\Gamma} = Fx\Gamma$. Moreover, $F/F \cap \Gamma$ has a finite invariant measure.*

Note that in the above theorem, the definition of a lattice contains the condition that $G/\Gamma$ has a finite volume. The first statement was called Raghunathan's topological conjecture, and the second statement was proved by Ratner in conjunction with Raghunathan's topological conjecture. Raghunathan's topological conjecture was published by Dani [7] for one dimensional unipotent groups and was generalized to groups generated by unipotent elements by Margulis [19] (before Ratner's proof). The proof for the general case was given by Ratner in a series of papers [25], [26], [27], [28]. For these, there is an excellent survey article by Ratner [29].

Note that in the above theorem, if $G$ is an algebraic group over $\mathbb{Q}$ and $\Gamma$ is an arithmetic lattice, the group $F$ becomes an algebraic group defined over $\mathbb{Q}$. For this, the reader should see Proposition (3.2) [34, pp. 321–322]. it is also proved in Proposition (3.2) [34, pp. 321–322] that the radical of $F$ is a unipotent subgroup. In [34], only one lattice is considered, but one can deduce the above statement for any lattice commensurable with the lattice in [34] by a simple argument using Ratner's theorem.

If we assume the above theorem, it is relatively easy to deduce the following corollary, which will immediately imply Theorem (0.1).

**Corollary (0.3)** *Let $Q$ be a real quadratic form in three variables as in* (0.1). *Then* $\mathrm{SO}(Q)_{\mathbb{R}}\mathrm{SL}(3)_{\mathbb{Z}}$ *is dense in* $\mathrm{SL}(3)_{\mathbb{R}}$.

Note that the representation $V = \mathrm{Sym}^2 k^n$ of $G = \mathrm{GL}(n)$ is an example of what we call prehomogeneous vector spaces. We recall the definition of prehomogeneous vector spaces.

**Definition (0.4)** *Let $G$ be a connected reductive group, $V$ a representation of $G$, and $\chi$ a non-trivial character of $G$, all defined over $k$. Then $(G, V, \chi)$ is called a prehomogeneous vector space if it satisfies the following properties.*
(1) *There exists a Zariski open orbit.*
(2) *There exists a polynomial $\Delta(x) \in k[V]$ such that $\Delta(gx) = \chi(g)^a \Delta(x)$ for a certain positive integer $a$.*

Such $\Delta$ is called a relative invariant polynomial. We define $V^{\mathrm{ss}} = \{x \in V \mid \Delta(x) \neq 0\}$ and call it the set of semi-stable points. If $(G, V, \chi)$ is an irreducible representation, the choice of $\chi$ is essentially unique and we may write $(G, V)$ as well. The theory of prehomogeneous vector spaces was initiated by Sato–Shintani [32] and Shintani [35]. If $(G, V)$ is irreducible, the classification is known (see [31]). For the zeta function theory of prehomogeneous vector spaces, see [36].

Now we can generalize the Oppenheim conjecture from our viewpoint in the following manner. For any algebraic group $G$, we denote the connected component of 1 in Zariski topology by $G^0$. If $G$ is defined over a subfield of $\mathbb{R}$, we denote the connected component of 1 of $G_{\mathbb{R}}$ in classical topology by $G^0_{\mathbb{R}+}$. Let $(G, V)$ be an irreducible prehomogeneous vector space over $\mathbb{Q}$. Let $H = [G, G] \subset G$ be the derived group. The group $H$ is semi-simple. Let $\Gamma \subset H^0_{\mathbb{R}+}$ be a lattice. For $x \in V^{\mathrm{ss}}_{\mathbb{R}}$, let $G_x$ be the stabilizer and $H_x = G_x \cap H$.

**Question (0.5)** (1) *If $x \in V^{\mathrm{ss}}_{\mathbb{R}}$ is sufficiently irrational, is $H^0_{x\mathbb{R}+}\Gamma$ dense in $H^0_{\mathbb{R}+}$?*
(2) *If* (1) *is true, find an explicit irrationality condition for* (1) *to be true.*

The Oppenheim conjecture can be considered as a special case of the above question for $G = \mathrm{GL}(n)$, $V = \mathrm{Sym}^2 \mathbb{Q}^n$. Suppose $(G_1, V_1)$, $(G_2, V_2)$ are prehomogeneous vector spaces, $\phi: G_1 \to G_2$ is a surjective homomorphism, and $\psi: V_1 \to V_2$ is a linear map, equivariant with respect to $\phi$. If $G_{2x}$ is finite for $x \in V^{\mathrm{ss}}_2$, Question (0.5)(1) is obviously false for the prehomogeneous vector space $(G_1, V_1)$. This eliminates cases (2) $n = 2$, (4), (8), (9), (11), (12), (28) in the classification (see [31, pp. 144–147]). Let $G, V, H$ be as before. We are interested in applying Ratner's theorem and if the stabilizer $H_x$ contains a split torus in its center, we cannot apply Ratner's theorem. This applies to cases (15) $m = 2$, (17), (18), (26). This leaves us twenty cases and one of them is case (2) $n \geq 3$, which led to the Oppenheim conjecture.

For these twenty cases with the split $k$–form of $G$, $H$ is either simple or a product of two simple groups. Since we are considering only algebraic groups, the center of each simple factor is finite. Consider $x \in V^{\mathrm{ss}}_{\mathbb{R}}$ such that the projection of $H^0_{x\mathbb{R}}$ to each simple factor of $H_{\mathbb{R}}$ has a positive rank. Note that if $H^0_{x\mathbb{R}}$ is split, this condition is satisfied (by case by case analysis). Then, by the Moore ergodicity theorem (see

Theorem 1 [20, p. 161]), the action of $H^0_{x\mathbb{R}_+}$ on $H^0_{\mathbb{R}_+}/\Gamma$ is ergodic (see (E) [20, p. 156] or [37, p. 8] for the definition of ergodicity). Therefore, by Proposition 2.1.7 [37, p. 10] for almost all $h\Gamma \in H^0_{\mathbb{R}_+}/\Gamma$ (meaning except for a measure zero set), $H^0_{x\mathbb{R}_+}h\Gamma$ is dense in $H^0_{\mathbb{R}_+}/\Gamma$. This implies that $H^0_{h^{-1}x\mathbb{R}_+}\Gamma$ is dense in $H^0_{\mathbb{R}_+}/\Gamma$. So starting from $x \in V^{ss}_{\mathbb{R}}$ as above, the answer to Question (0.5)(1) is affirmative for almost all points in the $H^0_{\mathbb{R}_+}$–orbit of $x$.

So in this series of papers, we consider Question (0.5)(2). In this paper, we consider the prehomogeneous vector space $G = \text{GL}(1) \times \text{GL}(8)$, $V = \wedge^3 k^8$. In part two, we consider prehomogeneous vector spaces $G = \text{GL}(1) \times \text{GL}(n)$, $V = \wedge^3 k^n$ for $n = 6, 7$, and $G = \text{GL}(2n)$, $V = \wedge^2 k^{2n}$.

Let $H_1 = \begin{pmatrix} 1 & & \\ & 1 & \\ & & -1 \end{pmatrix}$. We define

$$(0.6) \qquad W_{0\mathbb{R}} = \text{sl}(3)_{\mathbb{R}} = \{X \in \text{M}(3,3)_{\mathbb{R}} \mid \text{tr}(X) = 0\},$$
$$W_{1\mathbb{R}} = \text{su}(H_1)_{\mathbb{R}} = \{X \in \text{sl}(3)_{\mathbb{C}} \mid H_1{}^t\overline{X} + XH_1 = 0\}.$$

These are Lie algebras defined over $\mathbb{Q}$. We also define

$$(0.7) \qquad Q_0(X) = \text{tr}(X^2),\ F_0(X) = \det(X) \text{ for } X \in W_{0\mathbb{R}},$$
$$Q_1(X) = \text{tr}(X^2),\ F_1(X) = \sqrt{-1}\det(X) \text{ for } X \in W_{1\mathbb{R}}.$$

Obviously $Q_i$ is a quadratic form and $F_i$ is a cubic form on $W_{i\mathbb{R}}$. Our main result of this paper is the following

**Theorem (5.2)** *Let $i = 0$ or $1$ and $g \in \text{GL}(W_i)_{\mathbb{R}}$. Then if $Q_i(gX)$ is an irrational quadratic form in the sense of (0.1), the set of values of $F_i(gX)$ at primitive integer points is dense in $\mathbb{R}$.*

**Acknowledgement** The author would like to thank D. Witte for helpful conversations on ergodic theory, and R. Zierau for helping the author on representation theory. The proof of Proposition (4.3) is due to him.

## §1 Invariant theory of the space $\wedge^3 k^8$

Let $G = \text{GL}(1) \times \text{GL}(8)$, $V = \wedge^3 k^8$, and $\widetilde{T} = \text{Ker}(G \to \text{GL}(V))$. It is easy to see that
$$\widetilde{T} = \{(t^{-3}, tI_8) \mid t \in \text{GL}(1)\} \cong \text{GL}(1).$$

This is the prehomogeneous vector space of type (7) in [31] also discussed in Igusa [17].

Even though this is known to be a prehomogeneous vector space, we have to study the invariant theory of this space for our purposes.

Let $\{e_1, \cdots, e_8\}$ be a basis of $k^8$. We use the notation $e_{i_1 \cdots i_k} = e_{i_1} \wedge \cdots \wedge e_{i_k}$. It is known [31, p. 88] that the orbit of the element

$$(1.1) \qquad w = e_{123} + e_{456} + e_7 \wedge (e_{14} - e_{25}) + e_8 \wedge (e_{14} - e_{36})$$

is Zariski open in $V$.

Let

(1.2) $$\tau = \begin{pmatrix} 0 & I_3 & 0 \\ I_3 & 0 & 0 \\ 0 & 0 & -I_2 \end{pmatrix}.$$

Igusa [17, p. 275] pointed out that there is a split exact sequence

(1.3) $$1 \to G_w^0/\widetilde{T} \to G_w/\widetilde{T} \to \mathbb{Z}/2\mathbb{Z} \to 1,$$

where $\tau$ maps to the non-trivial element of $\mathbb{Z}/2\mathbb{Z}$ and

(1.4) $$G_w^0/\widetilde{T} \cong \mathrm{SL}(3)/Z(\mathrm{SL}(3)) \cong \mathrm{GL}(3)/Z(\mathrm{GL}(3)),$$
$$G_w/\widetilde{T} \cong \mathrm{Aut}\ (\mathrm{SL}(3)) \cong \mathrm{Aut}\ (\mathrm{GL}(3)),$$

where $Z(\mathrm{GL}(3)), Z(\mathrm{SL}(3))$ are the centers of $\mathrm{GL}(3), \mathrm{SL}(3)$ respectively. Here we interpret the quotient groups $\mathrm{SL}(3)/Z(\mathrm{SL}(3))$, $\mathrm{GL}(3)/Z(\mathrm{GL}(3))$ as those with universal categorical properties. In other words, if

$$\mathrm{SL}(3) = \mathrm{Spec}\ R_1,\ \mathrm{GL}(3) = \mathrm{Spec}\ R_2,\ H_1 = Z(\mathrm{SL}(3)),\ H_2 = Z(\mathrm{GL}(3)),$$

the ring of invariants $R_1^{H_1}$, $R_2^{H_2}$ are isomorphic. In this sense, $\mathrm{SL}(3)/Z(\mathrm{SL}(3)) \cong \mathrm{GL}(3)/Z(\mathrm{GL}(3))$. The point is that $k$–rational points of this quotient come from $\mathrm{GL}(3)_k$ but may not come from $\mathrm{SL}(3)_k$.

For any algebraic group $G$ over $k$, let $\mathrm{H}^1(k, G)$ be the first Galois cohomology set. We choose the definition so that trivial classes are those of the form $\{g^{-1}g^\sigma\}_{\sigma \in \mathrm{Gal}(\bar{k}/k)}$ ($g \in G_{\bar{k}}$) and the cocycle condition is $h_{\sigma\tau} = h_\tau h_\sigma^\tau$ for a continuous map $\{h_\sigma\}_{\sigma \in \mathrm{Gal}(\bar{k}/k)}$ from $\mathrm{Gal}(\bar{k}/k)$ to $G_{\bar{k}}$.

Note that by a similar argument as in [17], $G_w \cong \mathrm{GL}(3)/Z(\mathrm{GL}(3)) \times \widetilde{T}$. Since $\mathrm{H}^1(k, G) = \mathrm{H}^1(k, \widetilde{T}) = \{1\}$, the following theorem follows.

**Theorem (1.5) (Sato–Kimura, Igusa)**

$$G_k \setminus V_k^{\mathrm{ss}} \cong \mathrm{H}^1(k, \mathrm{Aut}\ (\mathrm{SL}(3))) \cong \mathrm{H}^1(k, \mathrm{Aut}\ (\mathrm{GL}(3))).$$

By the above theorem, the orbit space $G_k \setminus V_k^{\mathrm{ss}}$ corresponds bijectively with the set of $k$–isomorphism classes of $k$–forms of $\mathrm{GL}(3)$, $\mathrm{SL}(3)$. The correspondence is that if $x \in V_k^{\mathrm{ss}}$, $G_x^0/\widetilde{T}$ is isomorphic to the quotient of the corresponding $k$–form of $\mathrm{GL}(3)$ modulo its center.

The purpose of this section is to prove that we can associate a quadratic form and a cubic form in eight variables to any $x \in V_k^{\mathrm{ss}}$. We first describe how we can identify $G_w^0/\widetilde{T}$ with $\mathrm{GL}(3)/Z(\mathrm{GL}(3))$.

Let $W = \mathfrak{sl}(3)$ as in the introduction. We regard $V = \wedge^3 W$. Let $E_{ij}$ be the matrix whose only non-zero entry is the $(i,j)$–entry and is 1. We define

(1.6) $$H_1 = \frac{1}{3}\begin{pmatrix} 2 & & \\ & -1 & \\ & & -1 \end{pmatrix},\ H_2 = \frac{1}{3}\begin{pmatrix} 1 & & \\ & 1 & \\ & & -2 \end{pmatrix},$$
$$e_1 = E_{12},\ e_2 = E_{23},\ e_3 = -E_{31},\ e_4 = -E_{21},\ e_5 = -E_{32},\ e_6 = E_{13},$$
$$e_7 = H_1 - 2H_2 = \begin{pmatrix} 0 & & \\ & -1 & \\ & & 1 \end{pmatrix},\ e_8 = H_1 + H_2 = \begin{pmatrix} 1 & & \\ & 0 & \\ & & -1 \end{pmatrix}.$$

We choose $(e_1, \cdots, e_8)$ as the basis for $W$.

Let $\alpha_i, \beta_i, \gamma_i \in k$ for $i = 1, 2, 3$, and $\alpha_1 + \alpha_2 + \alpha_3 = 0$. We put

$$\alpha = (\alpha_1, \alpha_2, \alpha_3), \ \beta = (\beta_1, \beta_2, \beta_3), \ \gamma = (\gamma_1, \gamma_2, \gamma_3).$$

We define $H = H(\alpha)$, $X = X(\alpha, \beta, \gamma)$, $A = A(\alpha, \beta, \gamma)$ by

(1.7)
$$H = \alpha_1 H_1 + \alpha_2 H_2,$$
$$X = H - \beta_1 E_{12} - \beta_2 E_{23} - \beta_3 E_{31} - \gamma_1 E_{21} - \gamma_2 E_{32} - \gamma_3 E_{13},$$
$$A = \begin{pmatrix} \alpha_1 & 0 & 0 & 0 & \gamma_3 & \gamma_2 & \beta_1 & \beta_1 \\ 0 & \alpha_2 & 0 & -\gamma_3 & 0 & -\gamma_1 & -2\beta_2 & \beta_2 \\ 0 & 0 & \alpha_3 & -\gamma_2 & \gamma_1 & 0 & -\beta_3 & 2\beta_3 \\ 0 & -\beta_3 & -\beta_2 & -\alpha_1 & 0 & 0 & \gamma_1 & \gamma_1 \\ \beta_3 & 0 & \beta_1 & 0 & -\alpha_2 & 0 & -2\gamma_2 & \gamma_2 \\ \beta_2 & -\beta_1 & 0 & 0 & 0 & -\alpha_3 & -\gamma_3 & 2\gamma_3 \\ \gamma_1 & -\gamma_2 & 0 & \beta_1 & -\beta_2 & 0 & 0 & 0 \\ \gamma_1 & 0 & \gamma_3 & \beta_1 & 0 & \beta_3 & 0 & 0 \end{pmatrix}.$$

Then easy computations show that

(1.8) $$(\mathrm{ad}(X)(e_1), \cdots, \mathrm{ad}(X)(e_8)) = (e_1, \cdots, e_8) A.$$

This matrix $A$ is the one given in [31, p. 89]. Note that there is a misprint in [31, p. 89] and the $(3,7)$, $(8,1)$–entries are supposed to be $-\beta_3$, $\gamma_1$. It is proved in [31, p. 89] that if we regard $w \in \wedge^3 W$,

$$G_w^0 \cap \mathrm{GL}(W) = \{\mathrm{Ad}\,(g) \mid g \in \mathrm{GL}(3)\}.$$

As in [31], we define a map $D_3 : \wedge^3 W \to \wedge^2 W \otimes W$ by

(1.9) $$D_3(v_1 \wedge v_2 \wedge v_3) = v_2 \wedge v_3 \otimes v_1 - v_1 \wedge v_3 \otimes v_2 + v_1 \wedge v_2 \otimes v_3$$

for $v_1, v_2, v_3 \in W$.

For $x \in V_k$, we define

(1.10) $$S_x = x(\wedge, \otimes) D_3(x)(\wedge, \otimes) D_3(x) \in \wedge^7 W \otimes W \otimes W,$$

where $(\wedge, \otimes)$ means the wedge product for the first factor and the tensor product for the second factor.

Note that $\wedge^7 W \otimes W \otimes W \cong W^* \otimes (W \otimes W)$. So using the basis $\{e_1, \cdots, e_8\}$, we can regard $S_x$ as an $8 \times 8$ matrix with entries in the space of linear functions on $W$. We identify $W$ with the space of column vectors. It is easy to see that

(1.11) $$S_{(t,g)x}(v) = t^3 \det gg S_x(g^{-1}v)^t g$$

for $v \in W$.

**Definition (1.12)** $P_x(v) = \det S_x(v)$.

Clearly, $P_x(v)$ is bi-homogeneous of degree $(24, 8)$. It is easy to see that

(1.13) $$P_{gx}(v) = t^{24}(\det g)^{10} P_x(g^{-1}v).$$

Let

(1.14) $$w' = e_{123} + e_{156} + e_{246} + e_7 \wedge (e_{14} - e_{25}) + e_8 \wedge (e_{14} - e_{36}).$$

This is the element $X'_0$ considered in [31, p. 90].

Let $(f_1, \cdots, f_8)$ be the dual basis of $(e_1, \cdots, e_8)$. We computed $S_w, S_{w'}$ using the software "MAPLE" [1]. For example, in order to compute $S_w$, we associate a differential form $w = dx_1 \wedge dx_2 \wedge dx_3 + \cdots$. In order to distinguish two $D_3(w)$'s, we use variables $x_1, \cdots, x_8$ and $y_1, \cdots, y_8$. Then we type in the input as follows.

```
> with(difforms);
> defform(x1=0,x2=0,x3=0,x4=0,x5=0,x6=0,x7=0,x8=0);
> defform(y1=0,y2=0,y3=0,y4=0,y5=0,y6=0,y7=0,y8=0);
> w:= &^(d(x1),d(x2),d(x3))+&^(d(x4),d(x5),d(x6))
  +d(x7)&^(d(x1)&^d(x4)-d(x2)&^d(x5))
  +d(x8)&^(d(x1)&^d(x4)-d(x3)&^d(x6));
> v:= x1*d(x2)&^d(x3)-x2*d(x1)&^d(x3)+x3*d(x1)&^d(x2)
  +x4*d(x5)&^d(x6)-x5*d(x4)&^d(x6)+x6*d(x4)&^d(x5)
  +x7*d(x1)&^d(x4)-x1*d(x7)&^d(x4)+x4*d(x7)&^d(x1)
  -x7*d(x2)&^d(x5)+x2*d(x7)&^d(x5)-x5*d(x7)&^d(x2)
  +x8*d(x1)&^d(x4)-x1*d(x8)&^d(x4)+x4*d(x8)&^d(x1)
  -x8*d(x3)&^d(x6)+x3*d(x8)&^d(x6)-x6*d(x8)&^d(x3);
> u:= y1*d(x2)&^d(x3)-y2*d(x1)&^d(x3)+y3*d(x1)&^d(x2)
  +y4*d(x5)&^d(x6)-y5*d(x4)&^d(x6)+y6*d(x4)&^d(x5)
  +y7*d(x1)&^d(x4)-y1*d(x7)&^d(x4)+y4*d(x7)&^d(x1)
  -y7*d(x2)&^d(x5)+y2*d(x7)&^d(x5)-y5*d(x7)&^d(x2)
  +y8*d(x1)&^d(x4)-y1*d(x8)&^d(x4)+y4*d(x8)&^d(x1)
  -y8*d(x3)&^d(x6)+y3*d(x8)&^d(x6)-y6*d(x8)&^d(x3);
> w&^v&^u;
```

If the output contains $-x_1 y_7 dx_2 \wedge \cdots \wedge dx_8$ for example, we put $-f_1$ to the $(1, 7)$–entry of $S_w$. By interpreting the output of the above MAPLE session, $S_w$ turns out to be the following matrix

$$\begin{pmatrix} 0 & 3f_6 & 3f_5 & 3f_7 - 3f_8 & 0 & 0 & -f_1 & f_1 \\ 3f_6 & 0 & 3f_4 & 0 & 3f_8 & 0 & -f_2 & -2f_2 \\ 3f_5 & 3f_4 & 0 & 0 & 0 & -3f_7 & 2f_3 & f_3 \\ 3f_7 - 3f_8 & 0 & 0 & 0 & -3f_3 & -3f_2 & -f_4 & f_4 \\ 0 & 3f_8 & 0 & -3f_3 & 0 & -3f_1 & -f_5 & -2f_5 \\ 0 & 0 & -3f_7 & -3f_2 & -3f_1 & 0 & 2f_6 & f_6 \\ -f_1 & -f_2 & 2f_3 & -f_4 & -f_5 & 2f_6 & -2f_7 & 2f_8 - 2f_7 \\ f_1 & -2f_2 & f_3 & f_4 & -2f_5 & f_6 & 2f_8 - 2f_7 & 2f_8 \end{pmatrix}.$$

Similarly, $S_{w'}$ turns out to be the following matrix

$$\begin{pmatrix} 6f_7 - 6f_8 & 3f_6 & 3f_5 & 0 & -3f_3 & -3f_2 & -f_1 & f_1 \\ 3f_6 & -6f_8 & 3f_4 & 3f_3 & 0 & 3f_1 & -f_2 & -2f_2 \\ 3f_5 & 3f_4 & 0 & 0 & 0 & 0 & 2f_3 & f_3 \\ 0 & 3f_3 & 0 & 0 & 0 & -3f_5 & -f_4 & f_4 \\ -3f_3 & 0 & 0 & 0 & 0 & 3f_4 & -f_5 & -2f_5 \\ -3f_2 & 3f_1 & 0 & -3f_5 & 3f_4 & -6f_7 & 2f_6 & f_6 \\ -f_1 & -f_2 & 2f_3 & -f_4 & -f_5 & 2f_6 & -2f_7 & 2f_8 - 2f_7 \\ f_1 & -2f_2 & f_3 & f_4 & -2f_5 & f_6 & 2f_8 - 2f_7 & 2f_8 \end{pmatrix}.$$

We express elements of $W$ as $v = v_1 e_1 + \cdots + v_8 e_8$ or $v = \begin{pmatrix} v_1 \\ \vdots \\ v_8 \end{pmatrix}$. Let

(1.15)
$$\begin{aligned} Q_w(v) &= 2(-v_1 v_4 - v_2 v_5 - v_3 v_6 + v_7^2 - v_7 v_8 + v_8^2), \\ F_w(v) &= -v_1 v_2 v_3 + v_1 v_4 v_7 - v_1 v_4 v_8 + v_2 v_5 v_8 \\ &\quad - v_3 v_6 v_7 + v_4 v_5 v_6 - v_7^2 v_8 + v_7 v_8^2, \\ Q_{w'}(v) &= 2(-v_3^2 - v_4^2 + v_5^2), \\ F_{w'}(v) &= -v_1 v_3 v_5 + v_2 v_3 v_4 + v_3^2 v_7 + v_4 v_5 v_6 \\ &\quad - v_4^2 v_7 + v_4^2 v_8 + v_5^2 v_8. \end{aligned}$$

We can easily compute $\det S_w(v)$, $\det S_{w'}(v)$ by another MAPLE session as follows.

```
> with(linalg);
> a:= matrix(8,8);
> a[1,1]:= 0; ..., etc.
> factor(det(a));
```

It turns out that

$$P_w(v) = -1458 Q_w(v) F_w(v)^2, \quad P_{w'}(v) = -1458 Q_{w'}(v) F_{w'}(v)^2.$$

**Proposition (1.16)** (1) *There exist polynomials $Q_x(v), F_x(v)$ over $k$ such $Q_x(v)$ is bi-homogeneous of degree $(10, 2)$, $F_x(v)$ is bi-homogeneous of degree $(7, 3)$, $P_x(v) = -1458 Q_x(v) F_x(v)$, and $Q_w(v), F_w(v)$ coincide with the ones in (1.15).*
(2) $Q_{(t,g)x}(v) = t^{10} (\det g)^4 Q_x(g^{-1} v)$, $F_{(t,g)x}(v) = t^7 (\det g)^3 F_x(g^{-1} v)$.

*Proof.* We first assume that $k$ is algebraically closed. Let

$$\begin{aligned} X &= \operatorname{Spec} k[V, W]/(P_x(v)), \\ Y &= \operatorname{Spec} k[W]/(P_w(v)). \end{aligned}$$

By (1.13), $X$ is $G$–invariant. Let $\pi : X \to V$ be the projection map. Let $U = Gw$. Then $U$ is a Zariski open subset.

Consider the following map

$$f : G/G_w \times Y \ni (g, v) \to (gx, gv) \in \pi^{-1}(U) \subset X.$$

Then $f$ is set theoretically bijective to $\pi^{-1}(U)$. Therefore, $X$ is reducible.

Let $P_x(v) = p(x)q(x, v)$ where $q(x, v)$ is not divisible by a non-constant polynomial in $x$. The above observation implies that $q(x, v)$ is reducible. Since $P_w(v)$ does not have a linear factor, $q(x, v)$ does not have a linear factor either. Since $P_w(v)$ is a product of a quadratic polynomial and two cubic polynomials, $q(x, v)$ cannot have a factor of degree four with respect to $v$. So it has a unique factor, say $Q_x(v)$, of degree two with respect to $v$. Since $P_x(v)$ is bi-homogeneous, so is $Q_x(v)$. Multiplying a constant if necessary, we may assume that $Q_w(v)$ coincides with the one in (1.15).

Let $q(x, v) = Q_x(v)r(x, v)$. Consider the projection map

$$\pi : X_1 = \text{Spec } k[V, W]/(r(x, v)) \to V.$$

Then the fiber over $U$ does not have a reduced point. So by the generic smoothness [14, p. 272], $X_1$ cannot be irreducible. By considering the factorization of $P_x(v)$, there exist two polynomials $r_1(x, v), r_2(x, v)$, cubic with respect to $v$, such that $r(x, v) = r_1(x, v)r_2(x, v)$. Since $F_w(v)$ is irreducible and the zero set of $r_1(x, v), r_2(x, v)$ are the same generically, $r_1(x, v)$ is irreducible and is a constant multiple of $r_2(x, v)$. Multiplying a constant if necessary, there exists a polynomial, say $F_x(v)$, of degree three with respect to $v$ such that $r(x, v) = F_x(v)^2$ and $F_w(v)$ coincides with the one in (1.15). The polynomial $F_x(v)$ is also bi-homogeneous.

The discriminant of $Q_x(v)$ with respect to $v$ is a relative invariant polynomial and this reproves the existence of a relative invariant polynomial. Let $a, b$ be the degrees of $Q_x(v), F_x(v)$ with respect to $x$. Since $Q_{(t,g)x}(v), F_{(t,g)x}(v)$ are the unique quadratic and cubic factors of $P_{(t,g)x}(v)$, there exist integers $c, d$ such that

$$Q_{(t,g)x}(v) = t^a(\det g)^c Q_x(g^{-1}v), \; F_{(t,g)x}(v) = t^b(\det g)^d F_x(g^{-1}v).$$

This implies that if $p(x)$ is not constant, it is a relative invariant polynomial. However, $P_{w'}(v)$ is a non-zero polynomial and the discriminant of $Q_{w'}(v)$ is zero. So $p(x)$ must be a constant. By considering the point $w$, this constant must be $-1458$. Obviously $a + 2b = 24$ and $c + 2d = 10$. By considering scalar matrices, $3a = 8c - 2$, $3b = 8d - 3$. This system of linear equations has a unique solution $a = 10, b = 7, c = 4, d = 3$.

Now we consider an arbitrary field $k$ of characteristic zero. If $\sigma \in \text{Gal}(\bar{k}/k)$, by the uniqueness of the quadratic and cubic factors of $P_x(v)$, $Q_x(v)^\sigma$ and $F_x(v)^\sigma$ must be constant multiples of $Q_x(v)$ and $F_x(v)$ respectively. By considering the point $w$, these constants must be 1. So $Q_x(v)$ and $F_x(v)$ are Galois group invariant and therefore are polynomials over $k$ (ch $k = 0$). This proves the proposition. □

If $x \in V^{\text{ss}}$, $Q_x$ is non-degenerate. We define

(1.17) $$Q_x(v, v') = \frac{1}{2}(Q_x(v + v') - Q_x(v) - Q_x(v')).$$

By this paring, we get a map $i_x : W \to W^*$. This induces a map

(1.18) $$j_x : \wedge^3 W \to \wedge^3 W^*.$$

We define

(1.19) $$\Phi_x = j_x(x) \in \wedge^3 W^*.$$

Let $f_{ijk} = f_i \wedge f_j \wedge f_k$, etc. Easy computations show that

(1.20) $$\Phi_w = -f_{123} - f_{456} - f_7 \wedge (f_{14} - 2f_{25} + f_{36}) - f_8 \wedge (f_{14} + f_{25} - 2f_{36}).$$

**Proposition (1.21)** *For $(t,g) \in G$, $j_{(t,g)x}(v) = t^{30}(\det g)^{12} g j_x(g^{-1} v)$. In particular, $\Phi_{(t,g)x} = t^{31}(\det g)^{12} g \Phi_x$.*

*Proof.* If $(t, g) \in G$,
$$Q_{(t,g)x}(v, v') = t^{10}(\det g)^4 Q_x(g^{-1}v, g^{-1}v')$$
for all $v, v' \in W$. So
$$i_{(t,g)x}(v)(v') = t^{10}(\det g)^4 i_x(g^{-1}v)(g^{-1}v').$$

Considering the induced map on $\wedge^3 W$, $j_{(t,g)x}(v) = t^{30}(\det g)^{12} g j_x(g^{-1}v)$. This proves the proposition. □

The following proposition is an intuitive interpretation of $Q_w, F_w$. The proof is left to the reader.

**Proposition (1.22)** *For $X, Y, Z \in W$,*
$$Q_w(X, Y) = \mathrm{tr}(XY),$$
$$F_w(X) = \det X.$$

## §2 The fixed point set of $H_x^0$

Throughout this section, we assume that $k$ is algebraically closed. For our prehomogeneous vector space $(G, V)$, $H = \mathrm{SL}(W)$. We consider the fixed point set of $H_x^0$ in this section.

**Proposition (2.1)** *If $x \in V_k^{ss}$ and $y \in V_k$ is fixed by $H_{xk}^0$, $y$ is a scalar multiple of $x$.*

*Proof.* We first consider the case $x = w$.

**Lemma (2.2)** *If $y \in V_k$ is fixed by $H_{wk}^0$, $y$ is a scalar multiple of $w$.*

*Proof.* As we pointed out in §1, $H_k^0 \cong \mathrm{SL}(3)_k / Z(\mathrm{SL}(3)_k)$ (since $k$ is algebraically closed). Also as a representation of $\mathrm{SL}(3)$, $W$ is the adjoint representation, which has the weight $\Lambda_1 + \Lambda_2$. So we only have to show that $\wedge^3 W$ contains the trivial representation of SL(3) precisely once. We can verify this by the software "LIE" [6] as follows.

```
> setdefault A2
> alt_tensor(3,[1,1])
```

The result is
$$\wedge^3 W = U_1 \oplus U_2 \oplus U_3 \oplus U_4$$
where $U_1$ is the trivial representation and $U_2, U_3, U_4$ are irreducible representations with highest weights $3\Lambda_2, \Lambda_1 + \Lambda_2, 3\Lambda_1$ respectively. This proves the lemma. $\square$

We now consider the general case. Suppose $x = g_x w$ for $g_x \in G_k$. Then $H^0_{xk} = g_x H^0_{wk} g_x^{-1}$. So $g_x^{-1} y$ is fixed by $H^0_{wk}$. This implies that $g_x^{-1} y$ is a scalar multiple of $w$. Therefore, $y$ is a scalar multiple of $g_x w = x$. $\square$

**Corollary (2.3)** *If $y \in \wedge^3 W_k^*$ is fixed by $H^0_{wk}$, $y$ is a scalar multiple of $\Phi_w$.*

*Proof.* Suppose $y \in \wedge^3 W^*$ is fixed by $H^0_{wk}$. Let $z = j_w^{-1}(y)$. Since $H^0_w \subset \mathrm{SL}(W)$,
$$y = j_w(z) = j_{gw}(z) = g j_w(g^{-1} z) = gy = g j_w(z)$$
for $g \in H^0_w$. So $z$ is fixed by $H^0_{wk}$. By Lemma (2.2), $z$ is a scalar multiple of $w$. Therefore, $y$ is a scalar multiple of $\Phi_w$. $\square$

## §3 Lie algebra structures on $W$

In this section we construct Lie algebra structures on $W$ depending on $x \in V_k^{\mathrm{ss}}$. For $X, Y, Z \in \mathrm{Hom}(W, W) = W^* \otimes W$, we define

(3.1)
$$B(X, Y) = \mathrm{tr}(XY),$$
$$C(X, Y, Z) = \mathrm{tr}(XYZ - ZYX),$$
$$[X, Y] = XY - YX.$$

$B$ is a symmetric bilinear form and $C$ is an alternating tri-linear form. It is known [31, p. 90] that our prehomogeneous vector space has a relative invariant $\Delta(x)$ of degree 16 and therefore $\Delta((t, g)x) = t^{16} (\det g)^6 \Delta(x)$ for $(t, g) \in G$. Since we are assuming $\mathrm{ch}\, k = 0$, we may assume that $\Delta \in k[V]$ and $\Delta(w) = 1$. For $x \in V_k^{\mathrm{ss}}$, let $\mathfrak{h}_x$ be the Lie algebra of $H^0_x$.

**Proposition (3.2)** *Let $x \in V_k^{\mathrm{ss}}$. Then there exists a unique $k$–linear isomorphism $l_x : W \to \mathfrak{h}_x$ satisfying the following conditions*
(1) $l_x(Av) = [A, l_x(v)]$ for all $A \in \mathfrak{h}_x$, $v \in W$,
(2) $B(l_x(v), l_x(v')) = 6 Q_x(v, v')$ for all $v, v' \in W$,
(3) $\Delta(x) C(l_x(v), l_x(v'), l_x(v'')) = 6 \Phi_x(v, v', v'')$ for all $v, v', v'' \in W$.

*Proof.* We first consider the case $x = w$. Let

(3.3)
$$(\alpha, \beta, \gamma) = (v_7 + v_8, -2v_7 + v_8, -v_1, -v_2, v_3, v_4, v_5, -v_6),$$
$$l_w(\sum_{i=1}^{8} v_i e_i) = A(\alpha, \beta, \gamma),$$

where $A(\alpha, \beta, \gamma)$ is as in (1.7). Then $A(\alpha, \beta, \gamma)$ is the matrix representation of $\mathrm{ad}(v)$ for $v = \sum_{i=1}^{8} v_i e_i$. Therefore, $B$ is the Killing form and (2) is well known. Since

$$C(l_w(v), l_w(v'), l_w(v'')) = \mathrm{tr}(\mathrm{ad}(v)\mathrm{ad}(v')\mathrm{ad}(v'') - \mathrm{ad}(v'')\mathrm{ad}(v')\mathrm{ad}(v)),$$

this is an element of $\wedge^3 W^*$, fixed by $H^0_{w\bar{k}}$. So $C$ is a scalar multiple of $\Phi_w$ by Corollary (2.3). Let $A_i = l_w(e_i)$ for $i = 1, 2, 3$. We used MAPLE again to compute $C(A_1, A_2, A_3)$ and the result was $C(A_1, A_2, A_3) = -6$. Since $\Phi_w(e_1, e_2, e_3) = -1$, $C(l_w(v), l_w(v'), l_w(v'')) = 6\Phi_w(v, v', v'')$ for all $v, v', v'' \in W$. So we get (3).

Let $[\ ,\ ]_w$ be the Lie bracket for $W = \mathrm{sl}(3)$. Since $l_w(v)$ is the matrix representation of $\mathrm{ad}(v)$ and $[\mathrm{ad}(v), \mathrm{ad}(v')] = \mathrm{ad}([v, v']_w)$ for $v, v' \in W$,

(3.4)
$$l_w(v)v' = [v, v']_w,$$
$$[l_w(v), l_w(v')] = l_w([v, v']_w) = l_w(l_w(v)v').$$

Since $l_w$ is an linear isomorphism, any $A \in \mathfrak{h}_w$ is of the form $A = l_w(v')$. So (1) follows from the second equality in (3.4).

Suppose $l : W \to \mathfrak{h}_w$ is another map satisfying (1)–(3). Then there exists an element $g \in \mathrm{GL}(W)$ such that $l(v) = l_w(gv)$. By (1), $l_w(gAv) = [A, l_w(gv)]$ for all $A \in \mathfrak{h}_w, v \in W$. Let $v' = l_w^{-1}(A)$. Then since $l(v)$ also satisfies (1), by (3.4),

$$l_w(g[v', v]_w) = l_w(gl_w(v')v) = [l_w(v'), l_w(gv)] = l_w([v', gv]_w).$$

This is equivalent to
$$g[v', v]_w = [v', gv]_w$$

for all $v, v' \in W$. So $g$ commutes with the adjoint representation, which is irreducible. So by Schur's lemma, $g$ is a scalar matrix and $l(v) = cl_w(v)$ for a non-zero constant $c$. Since both $l_w, l$ satisfy properties (2), (3), $c^2 = c^3 = 1$. Therefore, $c = 1$. This proves the uniqueness.

Now we consider arbitrary $x \in V_k^{\mathrm{ss}}$. We first assume that $k$ is algebraically closed. So we choose $h \in \mathrm{GL}(W)_k$ so that $x = hw$. Let $l_x : W \to \mathfrak{h}_x$ be any linear isomorphism. Let
$$l'_x(v) = (\det h)^{-2} h^{-1} l_x(hv) h.$$

We show that the conditions (1)–(3) are equivalent to the following conditions

(3.5) $\quad l'_x((h^{-1}Ah)h^{-1}v) = [h^{-1}Ah, l'_x(h^{-1}v)],$
$\quad B(l'_x(h^{-1}v), l_x(h^{-1}v')) = 6Q_w(h^{-1}v, h^{-1}v'),$
$\quad C(l'_x(h^{-1}v), l'_x(h^{-1}v'), l'_x(h^{-1}v'')) = 6\Phi_w(h^{-1}v, h^{-1}v', h^{-1}v'').$

The equivalence of the first equations is straightforward. By the definition of $l'_x$,

$$B(l_x(v), l_x(v')) = (\det h)^4 B(hl'_x(h^{-1}v)h^{-1}, hl'_x(h^{-1}v')h^{-1})$$
$$= (\det h)^4 B(l'_x(h^{-1}v), l'_x(h^{-1}v')),$$
$$C(l_x(v), l_x(v'), l_x(v'')) = (\det h)^6 C(hl'_x(h^{-1}v)h^{-1}, hl'_x(h^{-1}v')h^{-1}, hl'_x(h^{-1}v'')h^{-1})$$
$$= (\det h)^6 C(l'_x(h^{-1}v), l'_x(h^{-1}v'), l'_x(h^{-1}v'')).$$

By Propositions (1.16)(2), (1.21),
$$Q_x(v, v') = (\det h)^4 Q_w(h^{-1}v, h^{-1}v'),$$
$$\Phi_x(v, v', v'') = (\det h)^{12} \Phi_w(h^{-1}v, h^{-1}v', h^{-1}v'').$$
Since $\Delta(x) = (\det h)^6 \Delta(w) = (\det h)^6$, the equivalences of the second and the third equations follow.

Since $A \in \mathfrak{h}_x$ is equivalent to $h^{-1}Ah \in \mathfrak{h}_w$, $l_w$ is the unique element which satisfies (3.5). Therefore, such $l_x$ exists and
$$(3.6) \qquad l_x(v) = (\det h)^2 h l_w(h^{-1}v) h^{-1}.$$

Now we consider arbitrary $k$. By the previous step, there exists a $\bar{k}$–linear isomorphism $l_x : W_{\bar{k}} \to \mathfrak{h}_{x\bar{k}}$ satisfying (1)–(3). Since $\Delta(x)$ and forms $B, Q_x, C, \Phi_x$ are all defined over $k$, for any Galois group element $\sigma$, $l_x^\sigma$ satisfies (1)–(3) also. Because of the uniqueness, $l_x^\sigma = l_x$ for all $\sigma \in \mathrm{Gal}(\bar{k}/k)$. Therefore, $l_x$ is defined over $k$. □

**Corollary (3.7)** *Let $x \in V_k^{\mathrm{ss}}$. Then, $l_x(gv) = g l_x(v) g^{-1}$ for all $g \in H_x^0, v \in W$.*

*Proof.* Let $F = \{g \in H_x^0 \mid l_x(gv) = g l_w(v) g^{-1}$ for all $v \in W\}$. Then $F$ is an algebraic subgroup of $H_x^0$. The tangent space of $F$ at $g = 1$ is a subspace of $\mathfrak{h}_x$ defined by the equation (3.2)(1) for all $v$. So the tangent space of $F$ at $g = 1$ is $\mathfrak{h}_x$. Since ch $k = 0$, $F$ must be $H_x^0$. □

**Definition (3.8)** *Let $x = hw \in V_k^{\mathrm{ss}}$ where $h \in \mathrm{GL}(W)_{\bar{k}}$. For $v, v' \in W_{\bar{k}}$ and $A \in \mathfrak{h}_{w\bar{k}}$, we define*
(1) $[v, v']_x = l_x^{-1}([l_x(v), l_x(v')])$,
(2) $m_{x,h}(v) = (\det h)^{-2} h v$,
(3) $n_{x,h}(A) = h A h^{-1}$.

It is easy to see that $[v, v']_x$ defines a Lie algebra structure on $W$ defined over $k$ and $l_x$ induces a Lie algebra isomorphism from $W$ to $\mathfrak{h}_x$. By Proposition (3.3), $Q_x, F_x$ can be regarded as functions on $\mathfrak{h}_x$ and $Q_x$ is $\frac{1}{6}$ times the Killing form of $\mathfrak{h}_x$.

**Proposition (3.9)** (1) $m_{x,h}([v, v']_w) = [m_{x,h}(v), m_{x,h}(v')]_x$ *for all $v, v' \in W_{\bar{k}}$ and therefore is a Lie algebra isomorphism.*
(2) *The following diagram is commutative.*

$$\begin{array}{ccc} & m_{x,h} & \\ W_{\bar{k}} & \longrightarrow & W_{\bar{k}} \\ l_w \downarrow & & \downarrow l_x \\ \mathfrak{h}_{w\bar{k}} & \longrightarrow & \mathfrak{h}_{x\bar{k}} \\ & n_{x,h} & \end{array}$$

*Proof.*
$$\begin{aligned} l_x(m_{x,h}([v, v']_w)) &= l_x((\det h)^{-2} h [v, v']_w) \\ &= h l_w([v, v']_w) h^{-1} \\ &= h [l_w(v), l_w(v')] h^{-1} \\ &= [h l_w(v) h^{-1}, h l_w(v') h^{-1}] \\ &= [l_x(m_{x,h}(v)), l_x(m_{x,h}(v'))]. \end{aligned}$$

Note that we used (3.4) in the third step. This proves (1). The proof of statement (2) is straightforward. □

It is possible to use $m_{x,h}$ to define a Lie algebra structure on $W$ and we don't need $\det h$ for that purpose. However, in order to define a map $l_x$ compatible with the group action, we need the factor $(\det h)^2$. In Proposition (3.2), we used the relative invariant of degree 16. But we could have used the discriminant $\Delta_1(x)$ of $Q_x(v)$ also. In that case since it is a relative invariant of degree 80, we use the condition
$$\Delta_1(x)C(l_x(v),l_x(v'),l_x(v''))^5 = 6^5\Phi_x(v,v',v'')^5.$$
Since 2 and 15 are relatively prime, the proof still works.

For the rest of this section, we assume $k = \mathbb{R}$ and we describe the orbit space $G_\mathbb{R}\backslash V_\mathbb{R}^{ss} \cong \mathrm{GL}(W)_\mathbb{R}\backslash V_\mathbb{R}^{ss}$. Let $\overline{X}$, etc. be the complex conjugate and $\sigma \in \mathrm{Gal}(\mathbb{C}/\mathbb{R})$ the corresponding Galois group element, i.e. $X^\sigma = \overline{X}$.

Since $\mathrm{SL}(3)$ is a simply connected algebraic group,
$$\mathrm{Aut}\,\mathrm{SL}(3) \cong \mathrm{Aut}\,\mathrm{sl}(3).$$

So instead of $\mathbb{R}$–forms of $\mathrm{SL}(3)$, we consider $\mathbb{R}$–forms of $\mathrm{sl}(3)$. It is well known that $\mathbb{R}$–forms of $\mathrm{sl}(3)$ are isomorphic to $\mathrm{sl}(3)$ or $\mathrm{su}(H)$ where $H$ is either one of the following two Hermitian matrices
$$H_1 = \begin{pmatrix} 1 & & \\ & 1 & \\ & & -1 \end{pmatrix},\quad H_2 = \begin{pmatrix} 1 & & \\ & 1 & \\ & & 1 \end{pmatrix}.$$

We define an $\mathbb{R}$–involution of $\mathrm{sl}(3)_\mathbb{C}$ by $p_H(X) = -H\,{}^t\overline{X}H^{-1}$. Then

(3.10) $$\mathrm{su}(H) = \{X \in \mathrm{sl}(3)_\mathbb{C} \mid p_H(X) = X\}.$$

If $x = hw$ for $h \in G_\mathbb{C}$, the corresponding cohomology class $c_x = \mathrm{H}^1(\mathbb{R},\mathrm{Aut}\,\mathrm{sl}(3))$ consists of a single element $h^{-1}h^\sigma = h^{-1}\overline{h} \in \mathrm{Aut}\,\mathrm{sl}(3)_\mathbb{C}$. Therefore, the corresponding $\mathbb{R}$–form is

(3.11) $$\mathfrak{h}_x \cong \{X \in \mathrm{sl}(3)_\mathbb{C} \mid h^{-1}\overline{h}(\overline{X}) = X\}.$$

So if we can find $h$ such that $p_H(X) = h^{-1}\overline{h}(\overline{X})$, the $\mathbb{R}$-orbit of $x = hw$ corresponds to $\mathrm{su}(H)$.

Let $D_1 = \begin{pmatrix} 1 & & \\ & -1 & \\ & & -1 \end{pmatrix}$ and $D_2 = I_3$. For $i = 1,2$, we define

(3.12) $$q_{D_i} = \begin{pmatrix} & D_i & \\ D_i & & \\ & & -I_2 \end{pmatrix},\quad h_{D_i} = \begin{pmatrix} I_3 & D_i & \\ \sqrt{-1}I_3 & -\sqrt{-1}D_i & \\ & & \sqrt{-1}I_2 \end{pmatrix}.$$

Note that $q_{D_i}^2 = I_8$, $H_i^2 = I_3$.

It is easy to see that $q_{D_i} X = -H_i{}^t X H_i^{-1}$. In other words, $q_{D_i}$ is the matrix representation of the involution $X \to -H_i{}^t X H_i^{-1}$ (this is not $p_H$). Also $q_{D_i} = h_{D_i}^{-1} \overline{h_{D_i}}$ for $i = 1, 2$. Therefore, $G_\mathbb{R} \backslash V_\mathbb{R}^{ss}$ is represented by $w_0 = w, w_1 = h_{D_1} w, w_2 = h_{D_2} w$. Easy computations show that

(3.13) $w_1 = 2(e_{123} - e_{156} + e_{246} - e_{345} + e_7 \wedge (e_{14} + e_{25}) + e_8 \wedge (e_{14} + e_{36}))$,
$w_2 = 2(e_{123} - e_{156} + e_{246} - e_{345} + e_7 \wedge (e_{14} - e_{25}) + e_8 \wedge (e_{14} - e_{36}))$.

**Lemma (3.14)** *For $i = 1, 2$, $q_{D_i} l_w(\overline{X}) q_{D_i}^{-1} = l_w(p_{H_i}(x))$.*

*Proof.* If $Y \in W$,

$$\begin{aligned} q_{D_i} l_w(\overline{X}) q_{D_i}^{-1} Y &= H_i{}^t [\overline{X}, H_i{}^t Y H_i^{-1}]_w H_i^{-1} \\ &= H_i [H_i^{-1} Y H_i, {}^t \overline{X}]_w H_i^{-1} \\ &= [Y, H_i{}^t \overline{X} H_i^{-1}]_w \\ &= [p_{H_i}(X), Y]_w \\ &= l_w(p_{H_i}(X)) Y. \end{aligned}$$

Therefore, $q_{D_i} l_w(\overline{X}) q_{D_i}^{-1} = l_w(p_{H_i}(x))$. □

Since $\mathfrak{h}_{w_i \mathbb{C}} = h_{D_i} \mathfrak{h}_{w \mathbb{C}} h_{D_i}^{-1}$ for $i = 1, 2$, Elements of $\mathfrak{h}_{w_i \mathbb{C}}$ can be written as $h_{D_i} l_w(X) h_{D_i}^{-1}$ for $X \in \mathrm{sl}(3)_\mathbb{C}$. We consider when this element belongs to $\mathfrak{h}_{w_i \mathbb{R}}$.

**Proposition (3.15)** *For $i = 1, 2$,*

$$\{h_{D_i}^{-1} X \mid X \in W_\mathbb{R}\} = \mathrm{su}(H_i)_\mathbb{R} \subset \mathrm{sl}(3)_\mathbb{C},$$
$$\{h_{D_i}^{-1} X \mid X \in W_\mathbb{Q}\} = \mathrm{su}(H_i)_\mathbb{Q} \subset \mathrm{sl}(3)_{\mathbb{Q}(\sqrt{-1})}.$$

*Proof.* Since $\overline{h_{D_i}} = h_{D_i} q_{D_i}$,

$$\begin{aligned} \overline{h_{D_i} l_w(X) h_{D_i}^{-1}} &= h_{D_i} q_{D_i} l_w(\overline{X}) q_{D_i}^{-1} h_{D_i}^{-1} \\ &= h_{D_i} l_w(p_{H_i}(X)) h_{D_i}^{-1}. \end{aligned}$$

Note that there is a natural $\mathbb{Q}$–structure on $\mathrm{su}(H_i)$. Therefore,

(3.16) $\mathfrak{h}_{w_i \mathbb{R}} = \{h_{D_i} l_w(X) h_{D_i}^{-1} \mid X \in \mathrm{su}(H_i)_\mathbb{R}\}$,
$\mathfrak{h}_{w_i \mathbb{Q}} = \{h_{D_i} l_w(X) h_{D_i}^{-1} \mid X \in \mathrm{su}(H_i)_\mathbb{Q}\}$.

It is easy to see that $\det h_{D_i} = -8\sqrt{-1}$ and $(\det h_{D_i})^2 = -64 \in \mathbb{Q}$ for $i = 1, 2$. Since $l_x$ is defined over any ground field of characteristic zero, $l_x(W_\mathbb{R}) = \mathfrak{h}_{x\mathbb{R}}$, $l_x(W_\mathbb{Q}) = \mathfrak{h}_{x\mathbb{Q}}$. Since the diagram (3.9)(2) commutes,

$$\{h_{D_i}^{-1} X \mid X \in W_\mathbb{R}\} = m_{x, h_{D_i}}^{-1}(W_\mathbb{R}) = \mathrm{su}(H_i)_\mathbb{R} \subset \mathrm{sl}(3)_\mathbb{C},$$
$$\{h_{D_i}^{-1} X \mid X \in W_\mathbb{Q}\} = m_{x, h_{D_i}}^{-1}(W_\mathbb{Q}) = \mathrm{su}(H_i)_\mathbb{Q} \subset \mathrm{sl}(3)_{\mathbb{Q}(\sqrt{-1})}.$$

This proves the proposition. □

By Proposition (1.16)(3),

(3.17)
$$Q_{w_i}(X) = 2^{12} Q_w(h_{D_i}^{-1}X),$$
$$F_{w_i}(X) = 2^9\sqrt{-1} F_w(h_{D_i}^{-1}X).$$

So if we identify $W_{\mathbb{R}} \cong \mathrm{su}(H_i)_{\mathbb{R}}$,

(3.18) $\quad Q_{w_i}(X) = 2^{12}\mathrm{tr}(X^2),\ F_{w_i}(X) = 2^9\sqrt{-1}\det(X).$

Note that if $X \in \mathrm{su}(H_i)_{\mathbb{R}}$, $\mathrm{tr}(X^2) \in \mathbb{R}$ and $\det(X)$ is purely imaginary.

It is easy to see that

$$h_{D_i}^{-1} = \begin{pmatrix} \frac{1}{2}I_3 & -\frac{\sqrt{-1}}{2}I_3 & \\ \frac{1}{2}D_i & \frac{\sqrt{-1}}{2}D_i & \\ & & -\sqrt{-1}I_2 \end{pmatrix}$$

for $i = 1, 2$. Easy computations show that

(3.20) $\quad Q_{w_1}(v) = 2^{11}(-v_1^2 + v_2^2 + v_3^2 - v_4^2 + v_5^2 + v_6^2 - 4v_7^2 + 4v_7v_7 - 4v_8^2),$
$\quad Q_{w_2}(v) = 2^{11}(-v_1^2 - v_2^2 - v_3^2 - v_4^2 - v_5^2 - v_6^2 - 4v_7^2 + 4v_7v_7 - 4v_8^2).$

So if $x \in G_{\mathbb{R}}w_i$, the signature of $Q_x(v)$ is $(5,3)$ or $(3,5)$, $(4,4)$, $(0,8)$ or $(8,0)$ for $i = 0, 1, 2$ respectively. This should be the case because these are constant multiples of the Killing forms of $\mathrm{sl}(3), \mathrm{su}(H_1), \mathrm{su}(H_2)$ respectively.

## §4 Intermediate groups

Let $\mathfrak{h}_1 = \mathrm{sl}(3)$, $\mathfrak{h}_2 = \mathrm{sl}(\mathfrak{h}_1)$, $\mathfrak{h}_3 = \mathrm{so}(Q_w)$. In this section, we consider Lie subalgebras of $\mathfrak{h}_2$ containing $\mathfrak{h}_1$.

First we assume that the ground field $k$ is algebraically closed. For a simple Lie algebra $\mathfrak{f}$ of rank $n$, let $\Lambda_1, \cdots, \Lambda_n$ be the fundamental weights of $\mathfrak{f}$ as in [31, pp. 7–32]. Let $d(\Lambda)$ be the dimension of the irreducible representation with highest weight $\Lambda$. The adjoint representation of $\mathfrak{h}_1$ is the irreducible representation with highest weight $\Lambda_1 + \Lambda_2$. This implies that its dual is equivalent to itself. So as a representation of $\mathfrak{h}_1$, $\mathfrak{h}_2$ plus the trivial representation is equivalent to $\mathfrak{h}_1 \otimes \mathfrak{h}_1$. We used LIE again to find the irreducible decomposition of this representation as follows.

```
> setdefault A2
> tensor([1,1],[1,1])
```

The result is that after eliminating the trivial representation, $\mathfrak{h}_2$ decomposes into a sum of five representations as follows

(4.1) $\qquad \mathfrak{h}_2 = U_1 \oplus U_2 \oplus U_3 \oplus U_4 \oplus U_5.$

where $U_1 = \mathfrak{h}_1$ and the highest weights of $U_2, U_3, U_4, U_5$ are $\Lambda_1+\Lambda_2, 3\Lambda_1, 3\Lambda_2, 2\Lambda_1+2\Lambda_2$ respectively. We can also find the dimensions of the representations as follows.

> dim([2,2]), etc.

The result is the following.

(4.2) $$\dim U_2 = 8, \ \dim U_3 = \dim U_4 = 10, \ \dim U_5 = 27.$$

Since $\mathfrak{h}_1 \subset \mathfrak{h}_3$ and $\dim \mathfrak{h}_3 = 28$, $\mathfrak{h}_3 = U_1 \oplus U_3 \oplus U_4$.

**Proposition (4.3)** *Suppose $k$ is algebraically closed. Then if $\mathfrak{h}_1 \subset \mathfrak{f} \subset \mathfrak{h}_2$ is a Lie subalgebra, $\mathfrak{f} = \mathfrak{h}_1, \mathfrak{h}_2$ or $\mathfrak{h}_3$.*

*Proof.* Suppose $\mathfrak{f} \neq \mathfrak{h}_1, \mathfrak{h}_2$. Then possible dimensions for $\mathfrak{f}$ are

(4.4) $$16, 18, 26, 28, 35, 36, 43, 45, 53, 55.$$

For simple Lie algebras of type $A_n, B_n, C_n, D_n$, the smallest non-trivial representation is the standard representation. The following table is the list of the dimensions of simple Lie algebras $\mathfrak{f}$ and the dimensions of the smallest non-trivial representations $V$ (see [15, pp. 461–474], [31, pp. 7–32]).

Table (4.5)

| Type | $\dim \mathfrak{f}$ | $\dim V$ |
|---|---|---|
| $A_n$ | $(n+1)^2 - 1$ | $n+1$ |
| $B_n$ | $2n^2 + n$ | $2n+1$ |
| $C_n$ | $2n^2 + n$ | $2n$ |
| $D_n$ | $2n^2 - n$ | $2n$ |
| $E_6$ | 78 | 27 |
| $E_7$ | 133 | 56 |
| $E_8$ | 248 | 3875 |
| $F_4$ | 52 | 26 |
| $G_2$ | 14 | 7 |

Lie algebras $\mathfrak{f}$ of types $A_n$–$D_n$ such that $\dim \mathfrak{f} \leq 55$ are the following

$A_n$: $3(n=1), 8, 15, 24, 35, 48$,
$B_n$: $10(n=2), 21, 36, 55$,
$C_n$: $10(n=2), 21, 36, 55$,
$D_n$: $28(n=4), 45$.

Suppose $\mathfrak{f}$ is simple. Then the only possibilities are $A_3, A_5, B_5, D_4, D_5, F_4, G_2$. Since $\mathfrak{h}_1 \subset \mathfrak{f} \subset \mathfrak{h}_2$ and $\mathfrak{h}_1$ is an irreducible eight dimensional representation of $\mathfrak{h}_1$, $\mathfrak{f}$ has an irreducible eight dimensional representation. By the above table, this is not possible for $B_5, C_5, F_4$. If $\dim \mathfrak{f} = 28$, the only possibility is $\mathfrak{f} = U_1 \oplus U_3 \oplus U_4 = \mathfrak{h}_3$.

Suppose $\mathfrak{f}$ is of type $A_3$. If the highest weight of a representation is $\Lambda = \sum_\nu m_\nu \Lambda_\nu$ and $\sum_\nu m_\nu \geq 2$, $d(\Lambda)$ is at least

$$d(2\Lambda_1) = d(2\Lambda_3) = 10,\ d(\Lambda_1 + \Lambda_3) = 15,$$
$$d(\Lambda_1 + \Lambda_2) = d(\Lambda_2 + \Lambda_3) = d(2\Lambda_2) = 20$$

(we used LIE again). So this cannot happen. Since $d(\Lambda_1) = d(\Lambda_3) = 4$, $d(\Lambda_2) = 6$, $\mathfrak{f}$ does not have an irreducible eight dimensional representation. Similarly, $\mathfrak{f}$ cannot be of type $A_5$ or $G_2$.

Suppose $\mathfrak{f}$ is semi-simple but not simple. Let $\mathfrak{f} = \mathfrak{f}' \times \mathfrak{f}''$, $\mathfrak{f}'$ simple, and $\mathfrak{f}''$ semi-simple. Then the eight dimensional irreducible representation of $\mathfrak{f}$ is a tensor product of irreducible representations $V', V''$ of $\mathfrak{f}', \mathfrak{f}''$. Since $\mathfrak{f}'$ is simple, there is only one one dimensional representation, which is the trivial representation. If $V'$ is the trivial representation, $\mathfrak{f}'$ acts trivially on $\mathfrak{h}_1$, which is a contradiction. Therefore, any simple factor of $\mathfrak{f}$ has an irreducible representation of dimension either two or four. The only such possibilities are type $A_1$ dimension 2 or 4, type $A_3$ dimension 4, or type $C_2$ dimension 4. So the only possibilities for $\mathfrak{f}$ are

$$\mathrm{sl}(2) \times \mathrm{sl}(2) \times \mathrm{sl}(2),\ \mathrm{sl}(2) \times \mathrm{sl}(4),\ \mathrm{sl}(2) \times \mathrm{sp}(4).$$

Then $\dim \mathfrak{f} = 9, 18, 13$. The first and the third cases cannot happen. Suppose $\mathfrak{f} = \mathrm{sl}(2) \times \mathrm{sl}(4)$. Since $\mathrm{sl}(3)$ does not have a non-trivial two dimensional representation, $\mathrm{sl}(4)$ alone contains $\mathrm{sl}(3)$. Since $\dim \mathfrak{f}$ cannot be 15, this cannot happen.

Suppose $\mathfrak{f}$ is reductive. By Schur's lemma, the center of $\mathfrak{f}$ must act trivially (there are no scalar matrices in $\mathfrak{h}_2$ except for the zero matrix). But this is a contradiction. Therefore, the only proper reductive subalgebra of $\mathfrak{h}_2$ containing $\mathfrak{h}_1$ is $\mathfrak{h}_3$.

It is known ([16, p. 187]) that any maximal proper connected subgroup of a reductive group is either reductive or parabolic. Therefore, any maximal proper subalgebra of a reductive Lie algebra is either reductive or parabolic. Suppose $\mathfrak{h}_1 \subset \mathfrak{f} \subset \mathfrak{h}_3$. If $\mathfrak{f} \neq \mathfrak{h}_1, \mathfrak{h}_3$, we choose $\mathfrak{f}$ so that it is maximal in $\mathfrak{h}_3$. Since we already proved that there is no such reductive subalgebra, $\mathfrak{f}$ must be a parabolic subalgebra of $\mathfrak{h}_3$. Then the reductive part of $\mathfrak{f}$ must be $\mathfrak{h}_1$. This is a contradiction because the rank of the Cartan subalgebra of $\mathfrak{h}_3$ is four, but the rank of the Cartan subalgebra of $\mathfrak{h}_1$ is three. So there is no proper subalgebra between $\mathfrak{h}_1$ and $\mathfrak{h}_3$.

If $\mathfrak{f}$ is not contained in $\mathfrak{h}_3$, choose a maximal proper subalgebra $\mathfrak{f} \subset \mathfrak{f}' \subset \mathfrak{h}_2$. Then $\mathfrak{f}'$ must be a parabolic subalgebra of $\mathfrak{h}_2$. The reductive part of $\mathfrak{f}'$ contains $\mathfrak{h}_1$ and it must be either $\mathfrak{h}_1$ or $\mathfrak{h}_3$. Again this is a contradiction by considering the ranks of the Cartan subalgebras. $\square$

**Corollary (4.6)** *Let $x \in V_\mathbb{R}^{\mathrm{ss}}$. If $H^0_{x\mathbb{R}+} \subset F \subset \mathrm{SL}(W)_\mathbb{R}$ is a connected closed subgroup, $F = H^0_{x\mathbb{R}+}, \mathrm{SO}(Q_x)_\mathbb{R}$ or $\mathrm{SL}(W)_\mathbb{R}$.*

*Proof.* Let $\mathfrak{f}$ be the Lie algebra of $F$. Then

$$\mathfrak{h}_{x\mathbb{R}} \subset \mathfrak{f} \subset \mathrm{sl}(W)_\mathbb{R}.$$

So

$$\mathfrak{h}_{x\mathbb{C}} \subset \mathfrak{f}_\mathbb{C} \subset \mathrm{sl}(W)_\mathbb{C}$$

and $\mathfrak{f} = \mathfrak{f}_{\mathbb{C}} \cap \mathrm{sl}(W)_{\mathbb{R}}$.

We choose $g \in G_{\mathbb{C}}$ so that $x = gw$. So

$$g^{-1}\mathfrak{h}_{x\mathbb{C}}g = \mathfrak{h}_{w\mathbb{C}} \subset g^{-1}\mathfrak{f}_{\mathbb{C}}g^{-1} \subset \mathrm{sl}(W)_{\mathbb{C}}.$$

By Proposition (4.3), $g^{-1}\mathfrak{f}_{\mathbb{C}}g^{-1} = \mathfrak{h}_{w\mathbb{C}}, \mathrm{so}(Q_w)_{\mathbb{C}}$, or $\mathrm{sl}(W)_{\mathbb{C}}$. This implies that $\mathfrak{f}_{\mathbb{C}}$ is $\mathfrak{h}_{x\mathbb{C}}, \mathrm{so}(Q_x)_{\mathbb{C}}$, or $\mathrm{sl}(W)_{\mathbb{C}}$. So $\mathfrak{f} = \mathfrak{f}_{\mathbb{C}} \cap \mathrm{sl}(W)_{\mathbb{R}} = \mathfrak{h}_{x\mathbb{R}}, \mathrm{so}(Q_x)_{\mathbb{R}}$ or $\mathrm{sl}(W)_{\mathbb{R}}$.

Since the correspondence between connected closed subgroups and Lie subalgebras is one-to-one, this proves the corollary. $\square$

## §5 An analogue of the Oppenheim conjecture

We prove our main theorem in this section.

**Theorem (5.1)** *Let $\Gamma \subset \mathrm{SL}(W)_{\mathbb{R}}$ be an arithmetic lattice and $x \in V_{\mathbb{R}}^{\mathrm{ss}}$ Suppose $Q_x(v)$ is indefinite and irrational in the sense of (0.1). Then $H_{x\mathbb{R}+}^0\Gamma$ is dense in $\mathrm{SL}(W)_{\mathbb{R}}$.*

*Proof.* Since $Q_x(v)$ is indefinite, $G_{x\mathbb{R}}^0$ is not compact and is generated by unipotent elements. Therefore, by Ratner's theorem, there exists a connected closed subgroup $H_{x\mathbb{R}+}^0 \subset F \subset \mathrm{SL}(W)_{\mathbb{R}}$ such that $\overline{H_{x\mathbb{R}+}^0\Gamma} = F\Gamma$. As we pointed out in the introduction, $F$ is defined over $\mathbb{Q}$. By Corollary (4.6), $F = H_{x\mathbb{R}+}^0, \mathrm{SO}(Q_x)_{\mathbb{R}}$, or $\mathrm{SL}(W)_{\mathbb{R}}$. We prove that the first two cases cannot happen.

Suppose $F = H_{x\mathbb{R}+}^0$. Let $\sigma \in \mathrm{Aut}\,(\mathbb{C}/\mathbb{Q})$. Since $F$ is defined over $\mathbb{Q}$, $(H_{x\mathbb{C}}^0)^\sigma = H_{x^\sigma\mathbb{C}}^0 = H_{x\mathbb{C}}^0$. So $H_{x\mathbb{C}}^0$ fixes $x^\sigma$. By Proposition (2.1), $x^\sigma$ must be a scalar multiple of $x$. So $x$ defines a $\mathbb{Q}$–rational point in $\mathbb{P}(V)_{\mathbb{R}}$. However, by the map $x \to Q_x$, it corresponds to the class of $Q_x$ in $\mathbb{P}(\mathrm{Sym}^2 W^*)_{\mathbb{R}}$. Since this map is defined over $\mathbb{Q}$, it is a $\mathbb{Q}$–rational point, and this is a contradiction.

Suppose $F = \mathrm{SO}(Q_x)_{\mathbb{R}}$. Then for any $\sigma \in \mathrm{Aut}\,(\mathbb{C}/\mathbb{Q})$,

$$\mathrm{SO}(Q_x^\sigma)_{\mathbb{C}} = (\mathrm{SO}(Q_x)_{\mathbb{C}})^\sigma = \mathrm{SO}(Q_x)_{\mathbb{C}}.$$

So $\mathrm{SO}(Q_x)_{\mathbb{C}}$ fixes $Q_x^\sigma$. Then it is well known that $Q_x^\sigma$ is a scalar multiple of $Q_x$. So the class of $Q_x$ is a $\mathbb{Q}$–rational point of $\mathbb{P}(\mathrm{Sym}^2 W^*)$, which is a contradiction. $\square$

Consider $Q_i(X)$, $F_i(X)$ in the introduction for $i = 0, 1$.

**Theorem (5.2)** *Let $i = 0$ or $1$ and $g \in \mathrm{GL}(W_i)_{\mathbb{R}}$. Then if $Q_i(gx)$ is an irrational quadratic form in the sense of (0.1), the set of values of $F_i(gx)$ at primitive integer points is dense in $\mathbb{R}$.*

*Proof.* Let $x = g^{-1}w_i$ where $g \in G_{\mathbb{R}}$ and $i = 0$ or $1$. By Proposition (1.16)(2),

$$Q_x(v) = (\det g)^{-4} Q_{w_i}(gv), \ F_x(v) = (\det g)^{-3} F_{w_i}(gv).$$

So $Q_x(v)$ is irrational if and only if $Q_{w_i}(gv)$ is irrational. $\{F_x(v) \mid v \in W_{\mathbb{Z}}\}$ is dense in $\mathbb{R}$ if and only if $\{F_{w_i}(gv) \mid v \in W_{\mathbb{Z}}\}$ is dense in $\mathbb{R}$. Since the argument is similar, we only consider the case $i = 1$. Note that $Q_{w_1}(X) = 2^{12}Q_1(X)$, $F_{w_1}(X) = 2^9 F_1(X)$. However, these constants have no effect on the statement.

We use the identification $W_{\mathbb{Q}} \cong \operatorname{su}(H_1)_{\mathbb{Q}}$. Let $L \subset W_{\mathbb{Q}}$ be the lattice which corresponds to $\operatorname{su}(H_1)_{\mathbb{Z}}$ by this identification. We define

(5.3) $$\Gamma = \{h \in \operatorname{SL}(W)_{\mathbb{R}} \mid hL \subset L\}.$$

This is an arithmetic lattice in $\operatorname{SL}(W)_{\mathbb{R}}$.

**Lemma (5.4)** *For any real number $r$, there exists $h \in \operatorname{SL}(W)_{\mathbb{R}}$ and an primitive element $u \in L$ such that $F_{h^{-1}x}(u) = r$.*

*Proof.* We may assume $x = \lambda w_1$ for a certain non-zero real number $\lambda$. Then

$$F_{h^{-1}x}(v) = \lambda^7 F_{h^{-1}w_1}(v) = \lambda^7 F_{w_1}(hv).$$

So $F_{h^{-1}x}(u) = r$ is equivalent to $F_{w_1}(hu) = \lambda^{-7} r$. Replacing $r$ if necessary, we may assume $\lambda = 1$. Choose a $\mathbb{Z}$–basis $(u_1, \cdots, u_8)$ of $L$ so that $F_{w_1}(u_1) = s \neq 0$. Let $t = (rs^{-1})^{\frac{1}{3}}$. Then

$$(u_1, \cdots, u_8) \to (tu_1, t^{-\frac{1}{7}} u_2, \cdots, t^{-\frac{1}{7}} u_8)$$

induces an element $h$ of $\operatorname{SL}(W)_{\mathbb{R}}$. Clearly, $F_{w_1}(hu_1) = r$. $\square$

For $h$ in the above lemma, we choose $h_1 \in H^0_{x\mathbb{R}}, h_2 \in \Gamma$ so that $h_1 h_2$ is close to $h$. Then
$$F_{(h_1 h_2)^{-1} x}(u) = F_{h_2^{-1} h_1^{-1} x}(u) = F_{h_2^{-1} x}(u) = F_x(h_2(u))$$
is close to $F_x(hu)$. Since $h_2 \in \Gamma$, $h_2(u) \in L$ is a primitive point also. This proves Theorem (5.2). $\square$

**Remark (5.5)** Consider the situation in Theorem (5.2). A simple modification of the above proof shows a slightly stronger statement as follows.

Let $i = 0, 1$. For any $r_1, \cdots, r_7 \in \mathbb{R}$ and $\epsilon > 0$, There exists a $\mathbb{Z}$–basis $(u_1, \cdots, u_8)$ of $W_{i\mathbb{R}}$ such that $|F_i(gu_j) - r_j| < \epsilon$ for $j = 1, \cdots, 7$.

Akihiko Yukie
Oklahoma State University
Mathematics Department
401 Math Science
Stillwater OK 74078-1058 USA
yukie@math.okstate.edu
http://www.math.okstate.edu/~yukie